\documentclass[12pt,leqno]{article}
\usepackage{amsopn,amsmath,amssymb,amscd,latexsym}
\usepackage[all]{xy}
\newcommand{\Q}{{\mathbb{Q}}}
\newcommand{\R}{{\mathbb{R}}}
\newcommand{\Z}{{\mathbb{Z}}}
\newcommand{\ok}{\overline{k}}
\newcommand{\tor}{\mathrm{tor}}
\newcommand{\End}{\mathrm{End}\,}
\newcommand{\Hom}{\mathrm{Hom}}
\newcommand{\Pic}{\mathrm{Pic}}
\newcommand{\Lh}{{\mathcal L}}
\newcommand{\Ph}{{\mathcal P}}
\newcommand{\oQ}{\overline{\Q}}
\newcommand{\verk}{\mbox{\scriptsize $\,\circ\,$}}
\newcommand{\halb}{\frac{1}{2}}
\newtheorem{theorem}{Theorem}
\newtheorem{lemma}{Lemma}
\newtheorem{prop}{Proposition}
\newenvironment{rem}{\noindent {\bf Remark}}{}

\newenvironment{proof}{\noindent {\bf Proof}}{\mbox{}\hfill$\Box$}
\begin{document}

\title{Linear relations among the values of canonical heights from the existence of non-trivial endomorphisms}
\author{Niko Naumann}
\maketitle

\begin{abstract}

  We study the interplay between canonical heights and endomorphisms of an abelian variety $A$ over a number field $k$. In particular we show that whenever the ring of endomorphisms defined over $k$ is strictly larger than $\Z$ there will 
be $\Q$-linear relations among the values of a canonical height pairing evaluated
at a basis modulo torsion of $A(k)$.
\end{abstract}

\section{Introduction}\label{sec1}
Let $A$ be an abelian variety over a number field $k$. In \cite{Ne} N\'eron
constructed a canonical pairing
\[
A(\ok)\times \hat{A}(\ok)\longrightarrow\R.
\]
The choice of a polarization then determines a height pairing $<,>$ on $A(\ok)$.
As observed in \cite{Ta} the Rosati involution  of an endomorphism is the adjoint with respect
to this canonical height pairing of the endomorphism acting on rational points.
We show that this forces linear relations among the values of the height pairing.\\
For the precise formulation let us first recall a part of the Albert classification (c.f. \cite{Mu}, p. 201 for more details): If $A$ is $k$-simple 
then $D$:=End$_k^0(A)$:=End$_k(A)\otimes\Q$ is a skew field carrying an involution, viz. the Rosati involution. Then $A/k$ is said to be of type I if $D$ is a totally real
number field; the involution is trivial in this case. Types II and III 
comprise quaternion algebras over totally real number fields. Finally, $A/k$ is of type IV if the center of $D$ is a CM field; the restriction of
the involution to this field is then complex conjugation.\\
 Denote by $W \subseteq \R$ the $\Q$-span of $\langle A(k) , A(k) \rangle$. For $r := $rk$(A(k))$ we have 

\begin{equation}
  \label{eq:1}
  \dim_{\Q} (W) \le \frac{r}{2} (r + 1 )
\end{equation}
because $\langle , \rangle$ is symmetric. Taking into account 
endomorphisms of $A$ we can prove the following bound:

\begin{theorem}
  \label{theorem1}
Assume $A$ is $k$-simple. Then $r$ is divisible by  $[\End^0_k (A) : \Q ]$
and we have
\[
\mbox{\rm dim}_\Q(W)\leq\frac{r}{2}(\frac{r}{[\mbox{\rm End}_k^0(A):\Q]}+\alpha)
\]
where $\alpha = 1 , \halb , - \halb$ or $0$ according to whether $\End^0_k (A)$ is of type I, II, III or IV in the Albert classification.
\end{theorem}

\begin{rem}
  If $r \neq 0$ and $\End^0_{k} (A) \neq \Q$ then dim$_\Q(W)<\frac{r}{2}(r+1)$, i.e. if there are non-torsion points in $A (k)$ and non-trivial (i.e. $\not\in\Z$ ) endomorphisms of $A$ defined over $k$ then the bound in theorem \ref{theorem1} is strictly sharper than the a priori bound (\ref{eq:1}). So the values of the height pairing on a basis modulo torsion of $A (k)$ satisfy non-trivial $\Q$-linear relations (inside $\R$).
\end{rem}

As an example we prove:

\begin{theorem} \label{theorem2}
Let $A/\Q$ be an abelian 
surface with real multiplication and let $D$ be the discriminant of $\End^0_\Q (A)$. Then there is a basis $P_1,...,P_{2n}\in A(\Q)\otimes\Q$ such that for any canonical height $h$ we have:
\[ h(P_{n+i})=Dh(P_i)\quad\mbox{ for }i=1,...,n \]
\end{theorem}

Examples of such surfaces are provided by modular abelian surfaces, c.f. also
the discussion at the end of section \ref{sec3}.\\
Next we generalize to higher dimensional abelian varieties the well known fact that the canonical height $h$ on an elliptic curve satisfies
\begin{equation}
  \label{eq:2}
  h \verk \varphi = (\deg \varphi) h
\end{equation}
for any endomorphism $\varphi$. For this we will have to deal simultaneously with heights associated to possibly different line-bundles. We are interested
only in heights $h_{\cal L}$ afforded by {\em symmetric} line-bundles ${\cal L}$ and
introduce
\[ SL(A/k):=\{ {\cal L} \in Pic(A_k) : [-1]^*({\cal L})\simeq {\cal L} \} \otimes\Q. \]

\begin{theorem}
  \label{theorem3}
Assume $A$ is $k$-simple of dimension $g$. Let $\Lh_1 , \ldots , \Lh_s$ be symmetric line-bundles on $A$ constituting a $\Q$-basis of $SL(A/k)$. Then there are quadratic forms
\[
\alpha_{ij} : \End^0_k (A) \longrightarrow \Q \; ; \; i,j = 1 , \ldots , s
\]
over $\Q$ such that for all $P \in A (k) \otimes \Q$ and $\varphi \in \End^0_k (A)$ we have
\[
h_{\Lh_i} (\varphi (P)) = \sum_j \alpha_{ij} (\varphi) h_{\Lh_j} (P) \; .
\]
Finally:
\[
\det (\alpha_{ij} (\varphi)) = \deg (\varphi)^{s/g} \; .
\]
\end{theorem}

Here deg is the degree extended to a homogeneous polynomial function on 
$ \End^0_k (A)$, c.f. \cite{Mu}, 19, Thm. 2.

Section \ref{sec2} contains the proofs of the above results.  We give some examples and the proof of theorem \ref{theorem2} in the last section.
I would like to thank V. Talamanca for making a copy of \cite{Ta} available to me and 
J. Cremona for pointing out \cite{FLSSSW}. I would also like to thank 
the referee for making numerous remarks which led to a substantial
improvement of the exposition.\\

\section{Proofs}\label{sec2}

We recall some fundamental facts about heights on abelian varieties over number fields, c.f. \cite{La}, Chapter 5. There is a canonical homomorphism
\[
\Pic (A_{\ok}) \longrightarrow \{ \mbox{quadratic functions} \; A (\ok) \to \R \} \; , \; \Lh \mapsto h_{\Lh}
\]
which is natural in $A$. The function $h_\Lh$ is called the canonical height associated to $\Lh$. In particular, for the Poincar\'e bundle $\Ph$ on $A \times \hat{A}$ we have
\[
h_{\Ph} : A (\ok) \times \hat{A} (\ok) \longrightarrow \R \; .
\]
Here $\hat{A}$ is the dual abelian variety of $A$. Recall the natural homomorphism
\[
\Pic (A_{\ok}) \longrightarrow \Hom_{\ok} (A , \hat{A}) \; , \; \Lh \mapsto \varphi_{\Lh} \; 
\]
with $\varphi_{\Lh}(a)=t_a^*(\Lh)\otimes\Lh^{-1}$ ($t_a$ is translation by $a$), 
see \cite{Mu}, II.8. For $\Lh \in \Pic (A_{\ok})$ we denote by
\[
L_{\Lh} (x,y) := h_{\Lh} (x+y) - h_{\Lh} (x) - h_{\Lh} (y)\; ; \; x,y\in A(\ok)
\]
the bilinear form associated to $h_{\Lh}$.
We say that a line-bundle $\Lh$ on $A$ is {\em symmetric} if [-1]$^*(\Lh)\simeq\Lh$.

\begin{lemma}
  \label{lemma2}
  \begin{enumerate}
  \item $h_{\Ph}$ is bilinear.
  \item For $\Lh \in \Pic (A_{\ok}) : L_{\Lh} (x,y) =  - h_{\Ph} (x , \varphi_{\Lh} (y))$ and if $\Lh$ is symmetric then $h_{\Lh} (x) = - \halb h_{\Ph} (x, \varphi_{\Lh} (x))$.
  \item For $\lambda \in \End_{\ok} (A) : h_{\Ph} ( \lambda x , y) = h_{\Ph} (x , \hat{\lambda} y)$.
  \item For the restriction of $h_{\Ph}$ to $A (k) \times \hat{A} (k)$ the kernels on both sides are exactly the torsion subgroups.
  \end{enumerate}
\end{lemma}

\begin{proof}
  \begin{enumerate}
  \item \cite{La}, Ch. 5, Prop. 4.3.
  \item \cite{La}, Ch. 5, Thm. 4.5.
  \item $h_{\Ph} (\lambda x,y) = h_{\Ph} ((\lambda \times 1) (x,y)) = h_{(\lambda \times 1)^* (\Ph)} (x,y) = h_{(1 \times \hat{\lambda})^* (\Ph)} (x,y) = h_{\Ph} (x , \hat{\lambda} y)$.
  \item As $A = \Hat{\Hat{A}}$ (over $k!)$ it suffices to consider the left kernel. $A (k)^{\tor}$ is orthogonal to $\hat{A} (k)$ because $\R$ is torsion free. Let $P \in A (k)$ be orthogonal to $\hat{A} (k)$. Choose $\Lh$ a symmetric ample line-bundle on $A$, defined over $k$. Then $h_{\Lh} (P) = - \halb h_{\Ph} (P , \varphi_{\Lh} (P)) = 0$ because $\varphi_{\Lh} (P) \in \hat{A} (k)$ and so $P \in A (k)^{\tor}$ by \cite{La}, Ch. 5, Thm. 6.1.\nolinebreak 
  \end{enumerate}\nopagebreak
\end{proof}

Fix some symmetric ample line-bundle $\Lh$ on $A$, defined over $k$. Associated to this is the height pairing $L_\Lh$ and we denote by $\langle , \rangle$ the extension of $L_{\Lh} \, |_{A (k)}$ to $V := A (k) \otimes \Q$. Furthermore, $\Lh$ determines an involution $'$ on $\End^0_k (A) := \End_k (A) \otimes \Q$, the Rosati involution:
\[
\varphi' := \varphi^{-1}_{\Lh} \hat{\varphi} \varphi_{\Lh} \quad \mbox{for} \; \varphi \in \End^0_k (A) \; .
\]
We have a natural representation
\begin{equation}
  \label{eq:3}
  \End_k (A) \longrightarrow \End_{\Z} (A (k))
\end{equation}
and as $A (k)$ is finitely generated
\[
\End_{\Q} (A (k) \otimes \Q) \simeq \End_{\Z} (A (k)) \otimes \Q \; .
\]
So from (\ref{eq:3}) we get a natural representation
\[
\phi : \End^0_k (A) \longrightarrow \End_{\Q} (V) \; .
\]

\begin{lemma} \label{lemma1}
  For $\varphi \in \End^0_k (A)$ and $v,w \in V$ we have
\[
\langle \phi(\varphi) v , w \rangle = \langle v , \phi(\varphi') w \rangle \; .
\]
\end{lemma}

\begin{proof}

For $v,w \in V$ and $\varphi \in \End^0_k (A)$ we compute:
  \begin{eqnarray*}
    \langle \phi (\varphi) v ,w \rangle & = & L_{\Lh} (\varphi v,w) = -h_{\Ph} (\varphi v , \varphi_{\Lh} w) \quad \mbox{by Lemma \ref{lemma2}, 2)}\\
& = & -h_{\Ph} (v , \hat{\varphi} \varphi_{\Lh} w) \quad \mbox{by Lemma \ref{lemma2}, 3)} \\
& = & -h_{\Ph} (v , \varphi_{\Lh} \varphi' w) = L_{\Lh} (v , \varphi' w) = \langle v , \phi (\varphi') w \rangle \; .
  \end{eqnarray*}
\end{proof}

\begin{proof} (of theorem \ref{theorem1}): As $A$ is $k$-simple
$D$:=End$^0_k(A)$ is a skew field acting on $V=A(k)\otimes\Q$. So
we have $V\simeq D^n$ as left $D$-modules for some $n\geq0$, hence $[D:\Q]$
divides $r=$dim$_\Q(V)$. Here we consider $D$ as a left $D$-module by
multiplication, as usual. The height pairing corresponds to a
$\Q$-linear map
\begin{equation}\label{D1}
D^n\otimes_{\Q}D^n\longrightarrow \R.
\end{equation}
We consider $D$ also as a {\em right} $D$-module by $xy:=y'x$ $(x,y\in D)$, where
$'$ is the Rosati involution on $D$. Then the content of lemma \ref{lemma1}
is that (\ref{D1}) factors over a $\Q$-linear map
\begin{equation}\label{D2}
D^n\otimes_D D^n\longrightarrow\R.
\end{equation}
We can identify
\[
D^n\otimes_D D^n\stackrel{\simeq}{\longrightarrow} M_n(D)\quad , \quad
(x_i)_i\otimes(y_j)_j\mapsto (x_i'y_j)_{ij}.
\]
The subspace of $D^n\otimes_D D^n$ spanned by $v\otimes w-w\otimes v$ is then 
identified with

\begin{equation}\label{ab}
T:=\{(x_{ij})\in M_n(D) : x_{ij}'=-x_{ji}\}\subset M_n(D).
\end{equation}

Since the height pairing is symmetric, the map in (\ref{D2}) factors
over $M_n(D)/T$, hence
\[
\mbox{dim}_{\Q}(W)\leq\mbox{dim}_{\Q}(M_n(D)/T).
\]
In the computation of the dimension of $M_n(D)/T$ we write
$|\cdot |$ as short-hand for dim$_\Q(\cdot)$ and put $S:=\{x\in D : x=x'\}$
and $\eta:=|S|/|D|$. Recall that $r=n|D|$.\\
From (\ref{ab}) we see that
\[
|T|=(1+...+(n-1))|D|+n(|D|-|S|)=\frac{n(n-1)}{2}|D|+n(|D|-|S|),
\]
hence
\begin{eqnarray}
|M_n(D)/T|=n^2|D|-(\frac{n(n-1)}{2}|D|+n|D|(1-\eta))\nonumber\\
=\frac{n(n+1)}{2}|D|-n|D|(1-\eta)\nonumber\\
=\frac{r}{2}(n+1-2(1-\eta))=\frac{r}{2}(\frac{r}{|D|}+2\eta-1).\nonumber\\
\nonumber\end{eqnarray}
This proves the theorem with $\alpha:=2\eta-1$. From \cite{Mu}, 21 
we know $\eta=1,3/4,1/4,1/2$ for $A$ of type I, II, III or IV.
Accordingly $\alpha=1,1/2,-1/2$ or $0$, as claimed.
\end{proof}

Recall the notation 
\[
SL(A/k)=\{ {\cal L}\in Pic(A_k) : [-1]^*({\cal L})\simeq {\cal L}\}\otimes\Q.
\]
The map
\[
\phi_A: SL(A/k)\longrightarrow\{\mbox{quadratic forms:} A(k)\otimes\Q\rightarrow\R\}\quad , \quad {\cal L}\mapsto h_{\cal L}
\]
is a transformation of contravariant functors on the category of abelian varieties
up to isogeny over $k$.\\
Now we explain when different line-bundles give rise to the same height. Note that the answer
is not immediate because we look at heights restricted to $A(k)$ only and
not on all of $A(\bar{k})$.

\begin{prop}
  \label{prop1}

The map $\phi_A$ is injective if and only if every non-trivial
isogeny factor of $A$ over $k$ has a $k$-rational point of infinite order.
\end{prop}

\begin{proof}

Assume first that $A$ has a non-trivial factor $B$ with $B(k)\otimes\Q=0$.
There is a surjective homomorphism
\[
\pi:A\longrightarrow B
\]
defined over $k$. Fix ${\cal M}$ a non-trivial symmetric line-bundle on $B$
defined over $k$. Then ${\cal L}:=\pi^*({\cal M})$ is non-trivial but 
\[
\phi_A({\cal L})=h_{\cal L}=h_{\cal M}\circ\pi=0 \quad\mbox{ on}\quad A(k)\otimes\Q,
\]
so $\phi_A$ is not injective.\\
For the converse assume that every factor of $A$ has a rational point
of infinite order and there is a non-trivial ${\cal L}\in SL(A/k)$ with $h_{\cal L}=0$
on $A(k)$. We derive a contradiction as follows:
There is a simple sub-variety $i:B\hookrightarrow A$ such that 
${\cal M}:=i^*({\cal L})$ is not trivial but $h_{\cal M}=h_{\cal L}\circ
i$ is zero on $B(k)\otimes\Q\neq 0$. As $B$ is simple, $\varphi_{\cal M}:
B\rightarrow\hat{B}$ is an isogeny inducing an isomorphism 
$B(k)\otimes\Q\stackrel{\simeq}{\rightarrow}\hat{B}(k)\otimes\Q$.
But this map has to be zero: For $P,Q\in B(k)$ we compute, using lemma
\ref{lemma2}:
\[
0=L_{\cal M}(P,Q)=-h_{\cal P}(P,\varphi_{\cal M}(Q)),
\]
hence $\varphi_{\cal M}(Q)\in\hat{B}(k)^{tor}$ and $\varphi_{\cal M}(Q)=0$
in $\hat{B}(k)\otimes\Q$.

\end{proof}

In order to prove theorem \ref{theorem3} we want to exploit the relation 
$h_{\cal L}\circ\alpha=h_{\alpha^*({\cal L})}$. Now, the assignment
${\cal L}\mapsto\alpha^*({\cal L})$ does not give an honest
action of End$_k^0(A)$ on $SL(A/k)$ because it is not additive in ${\cal L}$
and to proceed further we identify $SL(A/k)$ with the subspace of
symmetric elements
\[
S:=\{\alpha\in\mbox{\rm End}_k^0(A) : \alpha'=\alpha\}\subset\mbox{\rm End}_k^0(A)
\]
as follows:

\begin{lemma}
  \label{lemma4}
Let $A$ be an abelian variety over the perfect field $k$ and choose an ample line-bundle on $A$, defined over $k$. Denote by $'$ and $\lambda$ the associated involution and polarization, respectively. Then
\[
\psi : SL(A/k) \longrightarrow \End^0_k (A) \; , \; \Lh \longmapsto \lambda^{-1} \verk \varphi_{\Lh}
\]
identifies $SL(A/k)$ with $S$.
\end{lemma}

\begin{proof}
  Observing that $SL(A/k)\simeq (NS(A_{\bar{k}})\otimes\Q)^{G_k}$ this is standard for $k=\bar{k}$, c.f. \cite{Mu}, p. 190. The general case follows because $\psi$ is $G_k$-linear and $'$ commutes with the action of $G_k$, the absolute
Galois group of $k$.
This short proof was pointed out to me by the referee.
\end{proof}

In terms of this identification the sought for expression of $\alpha^*({\cal L})$ is the following:

\begin{lemma}
  \label{lemma5}
Assumptions and notations being as in Lemma \ref{lemma4}, given $\alpha \in \End^0_k (A)$ and $\Lh \in SL(A/k)$ we have
\[
\psi (\alpha^* (\Lh)) = \alpha' \psi (\Lh) \alpha \quad \mbox{in} \; S \; .
\]
\end{lemma}

\begin{proof}
  \begin{eqnarray*}
    \psi (\alpha^* (\Lh)) = \lambda^{-1} \verk \varphi_{\alpha^* (\Lh)} & = & \lambda^{-1} \verk \hat{\alpha} \verk \varphi_{\Lh} \verk \alpha = \alpha' \verk \lambda^{-1} \verk \varphi_{\Lh} \verk \alpha \\
 & = & \alpha' \psi (\Lh) \alpha \; .
  \end{eqnarray*}
\end{proof}

\begin{proof} (of theorem \ref{theorem3})
In the notation of the theorem we have

\[
h_{\Lh_i} (\varphi (P)) = h_{\varphi^* (\Lh_i)} (P) = h_{\psi^{-1} (\varphi' \psi (\Lh_i) \varphi)} (P) \; ,
\]
where the second equality results from lemma \ref{lemma5}.
So we define maps $\alpha_{ij} : \End^0_k (A) \to \Q$, recalling that the $\psi(\Lh_i)$ form a basis of $S$, by:
\begin{equation}
  \label{eq:8}
  \varphi' \psi (\Lh_i) \varphi = \sum_j \alpha_{ij} (\varphi) \psi (\Lh_j) \; \; \mbox{ for all } \, \varphi\in \End_k^0(A) 
\end{equation}
to get $h_{\Lh_i} (\varphi (P)) = h_{\psi^{-1} (\sum_j \alpha_{ij} (\varphi) \psi (\Lh_j))} (P) = h_{\sum_j \alpha_{ij} (\varphi) \Lh_j} (P) = \sum_j \alpha_{ij} (\varphi) h_{\Lh_j} (P)$, as desired.

It is clear from (\ref{eq:8}) that the $\alpha_{ij}$ are quadratic forms over $\Q$. Consider
\[
N : \End^0_k (A) \longrightarrow \End_{\Q} (S) \; , \; \varphi \longmapsto (s \mapsto \varphi' s \varphi) \; .
\]
Then $\det \verk N$ is a norm-form on $\End^0_k (A)$ over $\Q$. As $\deg$ is another such form and $\End^0_k (A)$ is simple there is a rational number $k$ such that $\det \verk N = \deg^k$ (\cite{Mu}, p. 179). Evaluating this on $n \in \Z$ gives
\[
(\det \verk N) (n) = \det (n^2) = n^{2s} = \deg (n)^k = n^{2gk} \; ,
\]
hence $k = s/g$ and $\det (\alpha_{ij} (\varphi)) = \deg (\varphi)^{s/g}$.
\end{proof}

\begin{rem}
  Let $\Lh$ be symmetric. One would like to have a formula
\[
h_{\Lh} (\varphi P) = f (\varphi) h_{\Lh} (P) \; , \; \mbox{all} \; P \in A (k)
\]
for a suitable $f (\varphi) \in \Q$. The injectivity of $\phi_A$ in proposition \ref{prop1} and the above theorem show that this will hold exactly for those $\varphi$ which have $\alpha_{1j} (\varphi) = 0$ for $j \neq 1$ in the notation of the theorem (and $\Lh_1 := \Lh$). These $\varphi$ form the intersection of quadratic hyper-surfaces inside $\End^0_k (A)$ and the dimension of this intersection most crucially depends on $\dim_{\Q} (S)$. 

For $\dim_{\Q} (S) = 1$ the above specializes to
\[
h_{\Lh} (\varphi P) = \deg (\varphi)^{1/g} h_{\Lh} (P) \; .
\]
This covers in particular the case of elliptic curves. A glance at the Albert classification reveals that $\dim_{\Q} (S) = 1$ if and only if $\End^0_k (A)$ is $\Q$, an imaginary quadratic field or a definite quaternion algebra over $\Q$.
\end{rem}

\section{Examples}\label{sec3}

Assume $A$ is $k$-simple of dimension $g$ and $K=$End$_k^0(A)$ is a quadratic field of discriminant
$D$. For simplicity assume also rk$(A(k))=2$. Then in a suitable basis of
$A (k) \otimes \Q$ the matrix of the height pairing will be given by 

\[ \left( \begin{array}{cc}

\alpha & \beta \\
\beta & D\alpha \end{array} \right)
\]

if $D>0$ and by

\[ \left( \begin{array}{cc}

\alpha & 0 \\
0  & -D\alpha \end{array} \right)
\]
if $D<0$ (for some $\alpha,\beta\in\R)$.
This illustrates theorem \ref{theorem1} in this case because $K$ is of type I if $D>0$ and of type IV otherwise. Furthermore, in case $D<0$, we have
\[ h\circ\varphi=(\deg(\varphi))^{1/g} h\]
for any $\varphi\in K$ as already noted above because the involution is complex conjugation in this case, so dim$_\Q(S)=1$.
If, however, $D>0$ then $S=K$ and for the base $1,\sqrt{D}\in S$ the matrix 
$(\alpha_{ij}(\varphi))$ of theorem \ref{theorem3} becomes

\[ \left( \begin{array}{cc}

a^2+Db^2 & 2ab \\
2abD & a^2+Db^2 \end{array} \right)
\]
\nopagebreak for $\varphi=a+b\sqrt{D}$.
So we will have a ``transformation formula''
$h\circ\varphi=f(\varphi)h$ exactly for those $\varphi$ satisfying $ab=0$.
This illustrates theorem \ref{theorem3} and the last remark of section \ref{sec2}.\\

Finally, it might be tempting to apply theorem \ref{theorem1} to the
following question of Serre (\cite{Se}, 3.8.):\newline

What is the transcendence degree of the field generated by the values of
B (=height pairing on an elliptic curve $E$ over $\Q$) ? Can it be $<r(r+1)/2$,
where r=rank($E(\Q)$) ? \newline

Clearly, theorem \ref{theorem1} does not give any improvement over the obvious bound because End$_{\Q}^0(E)=\Q$
even if $E$ has complex multiplication as we now assume, say End$_{\oQ}^0(E)=K$.
It is known that, for $\Q\subset L$, we will have End$_{L}^0(E)=K$ if and only 
if $K\subset L$. So we might try to apply theorem \ref{theorem1} to $E/K$, 
hoping that $E(\Q)\subset E(K)$ has small co-rank. For $A:=$Res$^K_{\Q}(E\otimes K)$
 (Weil's restriction of scalars) we have $E(K)\simeq A(\Q)$ and $A$ is 
isogeneous over $\Q$ to $E\times E'$, where $E'$ is the quadratic twist of $E$
corresponding to $\Q\subset K$. Now, $E'$ itself is $\Q$-isogeneous to $E$ as
follows from \cite{Mi2}, Thm. 3, remark 2. So rk($E(K)$)=$2$ rk($E(\Q)$) and
the ``effect'' of applying theorem \ref{theorem1} to $E/K$ is exactly 
compensated by this increase of the rank.\newline

It is however immediate that higher dimensional abelian varieties over $\Q$ can 
provide examples where the analogue of Serre's question has a positive answer.
We formulate this for surfaces only, i.e. we prove theorem \ref{theorem2}:\\
Let $A/\Q$ be an abelian 
surface with real multiplication, i.e. End$_{\Q}^0(A)$ is a real 
quadratic field of discriminant $D$, say. Fix some $\alpha\in$End$_{\Q}^0(A)$ with $\alpha^2=D$. Then there is a
basis $\{P_1,...,P_n,\alpha P_1,...\alpha P_n\}$ of $A(\Q)\otimes\Q$ and as
the Rosati involution is trivial on End$_{\Q}^0(A)$ lemma \ref{lemma1} gives
$h(\alpha P_i)=<\alpha P_i,\alpha P_i>=<DP_i,P_i>=Dh(P_i)$.\\

In \cite{FLSSSW} there are a number of explicit examples of genus 2 curves over $\Q$
whose Jacobians meet the assumptions of theorem \ref{theorem2} with $n=1$, many
of which are even absolutely simple. Of course, the generators given in table $3$ of \cite{FLSSSW} are not designed to satisfy the conclusion of theorem \ref{theorem2}.

\vspace*{0.5cm}

\hspace*{\fill} \\

Niko Naumann\\
Mathematisches Institut der WWU M\"unster\\
Einsteinstr. 62\\
48149 M\"unster\\
Germany\\

e-mail: naumannn@math.uni-muenster.de

\end{document}